\documentclass[a4paper, oneside, 10pt]{article}
\usepackage[english]{babel}
\usepackage{amsfonts}
\usepackage{amssymb}
\usepackage{graphics}
\usepackage{amsrefs}
\usepackage{latexsym}
\begin{document}

\title{{\bf{\Large{Generalized Chinese restaurant construction \\of exchangeable Gibbs partitions and related results.}}}\footnote{{\it AMS (2000) subject classification}. Primary: 60G58. Secondary: 60G09.} \footnote{{Research partially supported by MUR research grant n. 2006/134525.}}}
\author{\textsc {Annalisa Cerquetti}\footnote{Corresponding author. Istituto di
Metodi Quantitativi, Viale Isonzo,
25, 20133 Milano, Italy.
E-mail: {\tt annalisa.cerquetti@unibocconi.it}}\\
  \it{\small Bocconi University, Milano, Italy }}
\newtheorem{teo}{Theorem}
\date{}
\maketitle{}

\begin{abstract}

By resorting to sequential constructions of exchangeable random partitions (Pitman, 2006), and exploiting some known facts about generalized Stirling numbers, we derive a generalized Chinese restaurant process construction of exchangeable Gibbs partitions of type $\alpha$ (Gnedin and Pitman, 2006). Our construction represents the natural theoretical probabilistic framework in which to embed some recent results about a Bayesian nonparametric treatment of estimation problems arising in genetic experiment under Gibbs, species sampling, models priors.\\

\end{abstract}
%\end{frontmatter}
%\includegraphics[scale=0.3]{Piran07_98} % stone.jpg <--!!

\section{Introduction}

In two recent papers (Lijoi et al. 2007, 2008) a Bayesian prior to posterior analysis for the  subclass of exchangeable partitions in Gibbs form of type $\alpha \in (0,1)$, first introduced in Pitman (2003), and largely studied in Gnedin and Pitman (2006), has been proposed for a nonparametric treatment of some inferential problems arising in genetic experiments. 

Here, we derive a generalized group sequential construction of exchangeable partitions in Gibbs form of type $\alpha$ to place this theory in its natural probabilistic framework and to provide new insights on the derivation of relevant results. Our construction, which relies on known results of Pitman (2003, 2006) and Gnedin and Pitman (2006), has potential applications for investigating additional distributional results for quantities of statistical interest when exploiting in a Bayesian nonparametric perspective the theory of exchangeable partitions. 

Notice that, while in Lijoi et al. (2007, 2008) the treatment is in terms of {\it factorial coefficients}, our treatment is in terms of {\it Stirling numbers}, which naturally arise when summing over spaces of partitions with fixed number of blocks (see e.g. Pitman, 2006, Ch. 1) and also admit a generalized version. Moreover many convolution, multiplicative and multinomial formulas hold for Stirling numbers that greatly improved and simplify the presentation. 

The paper is organized as follows. In Section 2 we recall some preliminaries and basic definitions on rising factorials, random partitions and generalized Stirling numbers. In Section 3 we focus on infinite exchangeable partitions and derive a group sequential version of the Chinese restaurant process for Gibbs partitions of type $\alpha$. In Section 4 we show how to embed some results in Lijoi et al.  (2007, 2008) in our setting.

\section {Preliminaries and basic definitions}

We start by recalling some known facts about rising factorials and Stirling numbers which we will largely exploit in the following. A comprehensive reference for the role of these numbers in the theory of exchangeable random partitions is Pitman (2006). For the sake of clarity we strictly adopt his notations. \\

For $n=0,1,2,\dots,$ and arbitrary real $x$ and $h$, let $(x)_{n\uparrow h}$ denote the $n$th factorial power of $x$ with increment $h$ (also called generalized {\it rising} factorial)
\begin{equation}
\label{factdef}
(x)_{n \uparrow h}:= x(x+h)\cdots(x+(n-1)h)=\prod_{i=0}^{n-1}(x+ih)=h^n(x/h)_{n \uparrow},
\end{equation}
where $(x)_{n \uparrow}$ stands for $(x)_{n\uparrow 1}$, and $(x)_{0 \uparrow h}=x^h$, for which the following multiplicative law holds 
\begin{equation}
\label{multiplicative}
(x)_{n+r \uparrow h}=(x)_{n\uparrow h} (x +n h)_{r \uparrow h}.
\end{equation}
From e.g. Normand (2004, cfr. eq. 2.41 and 2.45) a binomial formula also holds, namely
\begin{equation}
\label{bino}
(x+y)_{n \uparrow h}=\sum_{k=0}^n {n \choose k} (x)_{k \uparrow h} (y)_{n-k \uparrow h},
\end{equation}
as well as a generalized version of the multinomial theorem, i.e.
\begin{equation}
\label{multi}
(\sum_{j=1}^p z_j)_{n \uparrow h}= \sum_{n_j \geq 0, \sum n_j=n} \frac{n!}{n_1!\cdots n_p!} \prod_{j=1}^p (z_j)_{n_j \uparrow h}.
\end{equation}
%Notice that (\ref{multiplicative}) and (\ref{multi}) make  unnecessary  In fact,
Notice that for $m_j>0$, for every $j$, and $\sum_j m_j=m$, an application of (\ref{multiplicative}) yields
\begin{equation}
\label{miaaa}
(z_j)_{n_j +m_j -1 \uparrow}=(z_j)_{m_j-1_\uparrow} (z_j +m_j -1)_{n_j \uparrow}
\end{equation}
and by (\ref{multi})
$$
\sum_{n_j \geq 0, \sum n_j=n} \frac{n!}{n_1!\cdots n_p!} \prod_{j=1}^p (z_j)_{n_j +m_j -1 \uparrow}=\prod_{j=1}^p (z_j)_{m_j -1 \uparrow}(\sum_{j=1}^p (z_j +m_j -1))_{n \uparrow}=
$$
$$
=\prod_{j=1}^p (z_j)_{m_j -1 \uparrow}(m +\sum_{j=1}^p z_j -p)_{n \uparrow},
$$
which agrees with the result in Lemma 1. in Lijoi et al. (2008).\\\\
Now recall that a {\it partition} of the finite set $[n]=(1,\dots, n)$ into $k$ blocks is an {\it unordered} collection of non-empty disjoint sets $\{A_1,\dots, A_k\}$ whose union is $[n]$, where  the blocks $A_i$ are assumed to be listed in order of appearance, i.e. in the order of their least elements. Recall also that the sequence $(|A_1|,\dots, |A_k|)$ of the sizes of blocks, $(n_1, \dots, n_k)$, defines a {\it composition} of $n$, i.e. a sequence of positive integers with sum $n$ and call $\mathcal{P}_{[n]}^k$ the space of all partitions of $[n]$ with $k$ blocks.

From Pitman (2006, cfr. eq. (1.7)) we know that the number of ways to partition $[n]$ into $k$ blocks and assign each block a $W$ combinatorial structure is given by
\begin{equation}
\label{bell}
B_{n,k}(w_\bullet):=\sum_{\{A_1,\dots, A_k\}\in \mathcal{P}_{[n]}^k }\prod_{i=1}^k w_{|A_i|}
\end{equation}
%\end{document}
which is a polynomial in variables $w_\bullet=(w_1,\dots, w_{n-k+1})$, known as the {\it $(n,k)$th partial Bell polynomial}.
Now, for each unordered partition of $[n]$ into $k$ disjoint non empty blocks there are $k!$ different {\it ordered} partitions of $[n]$ into $k$ such blocks, and corresponding to each composition $(n_1,\dots, n_k)$ of $n$ with $k$ parts, there are
$
{n \choose {n_1, \dots, n_k}}=n!\prod_{i=1}^k \frac{1}{n_i!}
 $
different {\it ordered} partitions  $(A_1, \dots, A_k)$ of $[n]$ with $|A_i|=n_i$. Hence the definition of $B_{n,k}(w_\bullet)$ as a sum of products over $\mathcal{P}_{[n]}^k$, translates in terms of sum over {\it compositions} of $n$ into $k$ parts as
$$B_{n,k}(w_\bullet)=\frac{n!}{k!} \sum_{(n_1,\dots, n_k)} \prod_{i=1}^k \frac{w_{n_i}}{n_i!}.
$$
%he relevant case for our discussion is $w_i=(1- \alpha)$. 
In what follows we will largely exploit the notion of {\it generalized Stirling numbers}, (for a comprehensive treatment see Hsu and Shiue, 1998; see also Pitman, 2006, eq. 1.19). For arbitrary distinct reals $\eta$ and $\beta$, these are the connection coefficients $S_{n,k}^{\eta, \beta}$ defined by
$$
(x)_{n \downarrow \eta}= \sum_{k=0}^n S_{n,k}^{\eta, \beta} (x)_{k \downarrow \beta} 
$$
%or equivalently (controllare...)
%\begin{equation}
%(x)_{n \uparrow -\eta}= \sum_{k=0}^n S_{n,k}^{-\eta, -\beta} (x)_{k \uparrow -\beta},
%\end{equation}
and correspond to
$$
S_{n,k}^{\eta, \beta}=B_{n,k}((\beta -\eta)_{\bullet -1 \downarrow \eta}).
$$

\noindent For $\eta=-1$, $\beta=-\alpha$, and $\alpha \in (-\infty, 1)$, $S_{n,k} ^{-1, -\alpha}$ is defined by
\begin{equation}
\label{unoalpha}
(x)_{n \uparrow 1}=\sum_{k=0}^{n} S_{n,k} ^{-1, -\alpha} (x)_{k \uparrow \alpha}.
\end{equation}
For $w_{n_i}=(1 -\alpha)_{n_i-1\uparrow}$ and $\alpha \in [0,1)$, equation (\ref{bell}) yields
\begin{equation}
\label{bellalpha}
B_{n,k}((1-\alpha)_{\bullet-1 \uparrow})=\sum_{\{A_1,\dots, A_k\}\in \mathcal{P}_{[n]^k} }\prod_{i=1}^k (1-\alpha)_{n_i-1\uparrow}=\frac{n!}{k!}\sum_{(n_1,\dots,n_k)}\prod_{i=1}^k \frac{(1-\alpha)_{n_i-1\uparrow}}{n_i!}=S_{n,k}^{-1,-\alpha}.
\end{equation}
{\bf Remark 1.} As recalled in the Introduction, in Lijoi et al. (2007, 2008) the treatment is in term of {\it generalized factorial coefficients}, which are the connection coefficients  $\mathcal{C}(n,k;\alpha)$ defined by
$$
(\alpha y)_{n\uparrow 1}=\sum_{k=0}^{n} \mathcal{C}(n,k;\alpha)(y)_{k\uparrow 1},
$$
(cfr. Charalambides, 2005). 
%neverthless the quantity $\mathcal{C}_{n,k}^\alpha/\alpha^k$ largely appears in their results and derivations. 
From (\ref{factdef}) and (\ref{unoalpha}), 
%  $(x)_{k \uparrow \alpha}=\alpha^k \left(\frac x\alpha\right )_{k \uparrow 1}$
if $x=y \alpha$ then
$$
(y \alpha)_{n \uparrow 1}= \sum_{k=0}^{n} S_{n, k}^{-1, -\alpha}(y \alpha)_{k \uparrow \alpha}=\sum_{k=0}^n S_{n,k}^{-1, -\alpha} \alpha^k
(y)_{k \uparrow 1},
$$
hence 
\begin{equation}
\label{coeff}
S_{n,k}^{-1, -\alpha}=\frac{\mathcal{C}_{n,k}^\alpha}{\alpha^k}.
\end{equation}
The representation (37) in Lijoi et al. (2008), (Toscano, 1939), also holds for generalized Stirling numbers with the obvious changes (cfr. e.g. Pitman, 2006, eq. 3.19). Additionally, specializing formula (16) in Hsu and Shiue (1998), the following convolution relation holds, which defines {\it non-central} generalized Stirling numbers
\begin{equation}
\label{convo}
S_{n,k}^{-1, -\alpha, \gamma}= \sum_{s=k}^{n} {n \choose s} S_{s,k}^{-1, -\alpha} (-\gamma)_{n-s \uparrow 1},
\end{equation}
and by (\ref{coeff}), 
$$
\mathcal{C}_{n,k}^{\alpha, \gamma}=\alpha^k S_{n,k}^{-1, -\alpha, \gamma}= \sum_{s=k}^{n} {n \choose s} \mathcal{C}_{s,k}^{\alpha} (-\gamma)_{n-s \uparrow 1}.
$$ 
Hence the following variation of equation (38) in Lijoi et al. (2008) defines {\it non-central} generalized Stirling numbers as connection coefficients, 
\begin{equation}
\label{noncentralsti}
(y \alpha-\gamma)_{n \uparrow 1}= \sum_{k=0}^{n} S_{n, k}^{-1, -\alpha, \gamma} \alpha^k
(y)_{k \uparrow 1}=\sum_{k=0}^{n} S_{n, k}^{-1, -\alpha, \gamma}
(y\alpha)_{k \uparrow \alpha} .
\end{equation}

\section {Generalized sequential constructions of exchangeable partitions}
Random partitions are random objects that arise in many contexts, {\it exchangeable} random partitions arise e.g. by sampling from random, almost surely discrete, probability measures. First recall that given a law $Q$ on the space $\mathcal{P}_1^\downarrow$ of decreasing sequences of positive numbers with sum 1, and a law $H(\cdot)$ on a Polish space $(S, \mathcal {S})$, a {\it random discrete} probability measure (RDPM) $P$ on $\mathcal{S}$ may always be defined as $P(\cdot)=\sum_{i=1}^\infty P_i \delta_{X_i}(\cdot)$, for $X_i$ iid $\sim H(\cdot)$ and $(P_i) \sim Q$. 
From Kingman's theory of exchangeable random partitions (Kingman, 1978), sampling from $P$ induces a random partition $\Pi$ of the positive integers $\mathbb{N}$ by the exchangeable equivalence relation
$i \approx j \Leftrightarrow X_i=X_j$, that is to say two positive integers $i$ and $j$ belong to the same block of $\Pi$ if and only if $X_i=X_j$, where $X_i|P$ are iid $\sim P$. It follows that, for each restriction $\Pi_n=\{A_1,\dots, A_k\}$ of $\Pi$ to $[n]=\{1,\dots, n\}$, and for each $n=1,2,\dots$, 
$$Pr(\Pi_n=\{A_1,\dots, A_k\})=p(n_1,\dots, n_k),
$$ where, for $j=1,2,\dots, k$, $n_j=|A_j|\geq 1$ and $\sum_{j=1}^k n_j=n$, for some non-negative symmetric function $p$ of finite sequences of positive integers called the {\it exchangeable partition probability function} (EPPF) determined by $\Pi$ (see Pitman, 2006, for a comprehensive account on exchangeable random partitions and related stochastic processes).\\

The particularly tractable class of exchangeable random partitions of {\it Gibbs form of type $\alpha$} has been first introduced in Pitman (2003) and then studied in Gnedin and Pitman (2006). Results and explicit forms for many of its EPPFs have been recently obtained in Ho et al. (2007). Here we recall the basic notions. An exchangeable random partition $\Pi_n$ of the first $n$ positive integers, is said to be of {\it Gibbs form} if for some nonnegative weights $W=(W_j)$ and $V=(V_{n,k})$ the EPPF of $\Pi$ can be expressed in the product form  
\begin{eqnarray}
\label{gibbs1}
p(n_1,\dots, n_k)=V_{n,k} \prod_{j=1}^k W_{n_j}
\end{eqnarray}
for all $1 \leq k\leq n$, and all compositions $(n_1,\dots, n_k)$ of $n$. Gnedin and Pitman (2006) show that to define an {\it infinite} random partition of $\mathbb{N}$, i.e. a sequence $(\Pi_n)$ consistent as $n$ varies, satisfying
\begin{equation}
\label{EPPFcc}
p({\bf n})=p({n_1, \dots, n_k})=\sum_{j=1}^{k({\bf n})+1} p({ n_1,\dots, n_j+1,\dots, n_k})
\end{equation}
for all compositions of $n$, the weights $(W_j)$ must be of the form 
$W_{n_j}=(1-\alpha)_{{n_j-1}\uparrow}$
for $\alpha \in [-\infty, 1)$, (with $W_j=1$ for every $j$, for $\alpha=-\infty$), and the weights $(V_{n,k})$ must be the solution to the {\it backward} recursion  $V_{n,k}=(n-\alpha k)V_{n+1, k}+V_{n+1,k+1}$  with $V_{1,1}=1$. They also obtain the solutions, for each $\alpha<1$, identifying the extreme points of the infinite dimensional simplex of the possible weights $V$, and deriving corresponding families of extreme partitions, in terms of the laws of the corresponding ranked atoms $(P_i)$. Their fundamental result, already stated without proof in Pitman (2003, cfr. Th. 8), is the following:\\\\
{\bf Theorem 2}. [Gnedin and Pitman, 2006; Th. 12] {\it Each exchangeable Gibbs partition of a fixed type $\alpha \in [-\infty, 1)$, i.e. characterized by an EPPF of the form
\begin{equation}
\label{gibbsGP}
p(n_1,\dots, n_k)=V_{n,k} \prod_{j=1}^k (1-\alpha)_{n_j-1 \uparrow}
\end{equation}
is a unique probability mixture of extreme partitions of this type, which are}\\\\
{\it \begin{tabular}{ll}
i) for $\alpha \in [-\infty, 0)$ & $PD(\alpha, m|\alpha|)$ partitions with $m=0,1,\dots,\infty$,\\
ii) for $\alpha=0$ & $PD(0,\theta)$ partitions with $\theta \in [0,\infty)$,\\
iii) for $\alpha \in (0,1)$ &  $PK(\rho_\alpha|t)$ partitions with $t \in [0,\infty).$
\end{tabular}\\\\
}
Recall that $PD(\alpha, \theta)$ stands for the two-parameter Poisson-Dirichlet distribution (Pitman and Yor, 1997) and $PK(\rho_\alpha|t)$ for the conditional Poisson-Kingman distribution derived from the stable subordinator, (cfr. Pitman, 2003).\\\\
%and for all $k\geq 1$, where $p({\bf n}^{j+})=p(n_1,\dots,n_j+1,\dots,n_k)$. 
Now it is known (cfr. Pitman, 2006, Ch. 3), that the consistency condition ({\ref{EPPFcc}) allows to derive a sequential construction of exchangeable partitions, known as {\it Chinese restaurant process construction} (first devised by Dubins and Pitman in 1986), that, in its more general formulation, (cfr. Ishwaran and James, 2003) is usually introduced as follows:\\

Given an infinite EPPF, $p({\bf n})=p(n_1, \dots, n_k)$, assume that an {\it unlimited} number of customers arrives sequentially in a restaurant with an {\it unlimited} number of circular tables, each capable of sitting an {\it unlimited} number of customers. Let the first customer to arrive be seated at the first table. For $n \geq 1$, given $(n_1, \dots, n_k)$, the placement of the first $n$ customers at $k$ tables, the $n+1$th customer is:\\\\a) seated at the table $j$, for $1\leq j\leq k_n$, provided $p({\bf n})>0$, with probability 
\begin{equation}
\label{predchi}
p_{j,n}=p_j({\bf n})=\frac{p({\bf n}^{j+})}{p({\bf n})}
%=\frac{V_{n+1,k}}{V_{n,k}}(1-\alpha),
\end{equation}
where $p({\bf n}^{j+})$ stands for $p(n_1, \dots, n_j+1, \dots, n_k)$ or is\\\\
b) seated at a {\it new} table with probability 
\begin{equation}
\label{predchi2}
p_{0,n}=p_{0,n}({\bf n})=\frac{p({\bf n}^{l+})}{p({\bf n})}
%=\frac{V_{n+1,k+1}}{V_{n,k}}=1 -\frac{(n-\alpha k)V_{n+1,k}}{V_{n,k}}
\end{equation} 
for $l=k_n+1$, and $
\sum_{j=1}^{k+1} p_{j,n} + p_{0,n}=1.
$\\\\
Here we derive a {\it group} sequential version of  %rules (\ref{predchi}) and (\ref{predchi2}) 
the Chinese restaurant process which can be introduced as follows:\\

Given an infinite EPPF $p(\bf n)$, assume that an {unlimited} numbers of {\it groups} of customers arrive sequentially in a restaurant with an unlimited numbers of circular tables, each capable of sitting an unlimited numbers of customers.
Given the placement of the first group of $n$ in a $(n_1, \dots, n_k)$ configuration in $k$ tables, the {\it new} group of $m$ customers is: \\\\
a) seated at the {\it old} $k$ tables in configuration $(m_1,\dots, m_k)$, for $m_j \geq 0$, $\sum_{j=1}^k m_j=m$, with probability 
\begin{equation}
\label{allold}
p_{\bf m}({\bf n})=p({\bf m}| {\bf n})=\frac{p(n_1+m_1, \dots, n_k+m_k)}{p(n_1,\dots, n_k)},
\end{equation}
b) seated at $k^*$ {\it new} tables in configuration $(s_1, \dots, s_{k^*})$, for $\sum_{j=1}^{k^*} s_j =m$, $1 \leq k^* \leq m$, $s_j \geq 1$, with probability  
\begin{equation}
\label{allnew}
p_{\bf s}({\bf n})=p({\bf s}|{\bf n})\frac{p(n_1,\dots,n_k, s_1,\dots, s_{k^*})}{p(n_1,\dots, n_k)},
\end{equation}
c) $s < m$ customers are seated at $k^*$ {\it new} tables in configuration $(s_1,\dots,s_{k^*})$ and the remaining $m-s$ customers at the {\it old} tables in configuration $(m_1,\dots, m_k)$ for $\sum_{j=1}^{m} m_j= m-s$, $1 \leq s \leq m$, $\sum_{j=1}^{k^*} s_j=s$, $m_j \geq 0$, $s_j \geq 1$ with probability
\begin{equation}
\label{oldnew1}
p_{{\bf m, \bf s}}({\bf n})=p({\bf m}, {\bf s}| {\bf n})=\frac{p(n_1+m_1, \dots, n_k+m_k, s_1, \dots, s_{k^*})}{p(n_1,\dots, n_k)}.\\\\
\end{equation}

\noindent From now on we focus on the particular case of Gibbs EPPFs. The mathematical tractability of the Gibbs product form, combined with the properties of generalized rising factorials previously recalles, allows an easy derivation.\\\\
%to deduce from the previous definition the case of our interest.\\\\ 
{\bf Proposition 3.} For the Gibbs EPPF (\ref{gibbsGP}), formula (\ref{allold}), (\ref{allnew}) and (\ref{oldnew1}) specialize as follows.  Given the placement of the first group of $n$ customers in a $(n_1, \dots, n_k)$ configuration in $k$ tables, the new group of $m$ customers is \\\\
a) seated at the $k$ old tables in configuration $(m_1,\dots, m_k)$, for $m_j \geq 0$, $\sum_{j=1}^k m_j=m$, with probability 
\begin{equation}
p_{\bf m}({\bf n})=\frac{V_{n+m,k}\prod_{j=1}^k (1-\alpha)_{n_j+m_j-1}}{V_{n,k}\prod_{j=1}^k (1-\alpha)_{n_j -1\uparrow}}\nonumber
\end{equation}
which, by means of the multiplicative relation (\ref{miaaa}), simplifies to 
\begin{equation}
\label{gibbsallold}
=\frac{V_{n+m,k}}{V_{n,k}} \prod_{j=1}^k (n_j-\alpha)_{m_j\uparrow},
\end{equation}
b) seated at $k^*$ {\it new} tables with configuration $(s_1, \dots, s_{k^{*}})$, for $\sum_{j=1}^{k^*} s_j =m$, $1 \leq k^* \leq m$, $s_j \geq 1$, with probability  
\begin{equation}
\label{gibbsallnew}
p^{\bf s}({\bf n})=\frac{V_{n+m,k+k^*}\prod_{j=1}^k (1-\alpha)_{n_j -1\uparrow}\prod_{j=1}^{k^*} (1-\alpha)_{s_j -1\uparrow}}{V_{n,k} \prod_{j=1}^k (1-\alpha)_{n_j -1\uparrow}}=\frac{V_{n+m, k+k^*}}{V_{n,k}}\prod_{j=1}^{k^*} (1-\alpha)_{s_j -1\uparrow},\\
\end{equation}
c) a subset $s < m$ of the new customers is seated at $k^*$ {\it new} tables in configuration $(s_1,\dots,s_{k^*})$ and the remaining $m-s$ customers are seated at the {\it old} tables in configuration $(m_1,\dots, m_k)$ for $\sum_{j=1}^{k} m_j= m-s$, $1 \leq s \leq m$, $\sum_{j=1}^{k^*} s_j=s$, $m_j \geq 0$, $s_j \geq 1$ with probability
$$
%\label{gibbsoldnew}
p_{\bf m}^ {\bf s}({\bf n})=\frac{V_{n+m,k+k^*} \prod_{j=1}^k (1-\alpha)_{n_j+m_j-1\uparrow}\prod_{j=1}^{k^*}(1-\alpha)_{s_j-1\uparrow}}{V_{n,k}\prod_{j=1}^k (1-\alpha)_{n_j-1\uparrow}}=
$$
which simplifies to
\begin{equation}
\label{oldenew}
=\frac{V_{n+m,k+k^*}}{V_{n,k}}\prod_{j=1}^k (n_j-\alpha)_{m_j\uparrow}\prod_{j=1}^{k^*}(1-\alpha)_{s_j-1\uparrow}.\\
\end{equation}
{\bf Corollary 4.} Consider the event A=$\{$All $m$ new customers are seated at {\it new} tables$\}$, by (\ref{bellalpha}) and (\ref{gibbsallnew}), summing over all the way to allocate $m$ in $k^*$ tables for every $k^*$ yields
%B=$\{$A subset of $s$, for $s \in (1,\dots,m)$ customers seat at new tables and the remaining $(m-s)$ seat to old tables$\}$, C=$\{$New customers all seat at old $k$ tables$\}$. \\\\The probability all in new, $P(C)$ by By equation (13)
\begin{equation}
\label{proballnew}
Pr(A|n_1,\dots, n_k)=\sum_{k^*=1}^m \frac{1}{k^*!}\sum_{s_j\geq 1, \sum_j s_j=m}{m \choose {s_1,\dots, s_{k^*}}} \frac{V_{m+n,k+k^*}}{V_{n,k}} \prod_{j=1}^{k^*} (1-\alpha)_{s_j-1 \uparrow}=
\end{equation}
$$
=\sum_{k^* =1}^m \frac{V_{n+m, k+k^*}}{V_{n,k}} S_{m,k^*}^{-1,-\alpha}.
$$
For the event B=$\{$All $m$ new customers are seated at the  $k$ {\it old} tables$\}$, by (\ref{gibbsallold}) it is easy to obtain
\begin{equation}
\label{proballold}
Pr(B|n_1,\dots,n_k)
=\frac{V_{n+m,k}}{V_{n,k}} \mathop{\sum_{(m_1,\dots, m_k)}}_{\sum_j m_j=m, m_j\geq 0}{m \choose {m_1,\dots, m_k}} \prod_{j=1}^k (n_j -\alpha)_{m_j \uparrow}=
\end{equation}
%\mathop{\sum_{(m_1,\dots, m_k)}}_{\sum_j =m_j, m_j\geq 0}{m \choose {m_1,\dots, m_k}} \frac{V_{n+m,k}}{V_{n,k}}\frac{\prod_{j=1}^k (1-\alpha)_{n_j+m_j-1\uparrow})}{\prod_{j=1}^k \prod_{j=1^k}(1-\alpha)_{n_j-1\uparrow}}=
%\end{equation}
 and specializing (\ref{multi}) for $z_j=(n_j-\alpha)$, 
$$= \frac{V_{n+m,k}}{V_{n,k}} (\sum_{j=1}^k (n_j -\alpha))_{m \uparrow}= \frac{V_{n+m,k}}{V_{n,k}} (n -k\alpha)_{m\uparrow}.
$$

\section {Embedding the Bayesian nonparametric analysis in the group sequential construction.}

Now we show how to embed in our construction some of the results in Lijoi et al. (2007, 2008). Notice that our derivation shows explicitily how generalized Stirling numbers arise in this context, hence indirectly clarifies how generalized factorial coefficients arise in their treatment. \\\\
{\bf Proposition 5.} {The result in Proposition 1. of Lijoi et al. (2008) may be obtained by (\ref {oldenew}) by summing over the ways to choose $(m-s)$ integers from $m$, and all the ways to allocate $(m-s)$ integers in the $k$ old tables}.\\\\
{\it Proof.} Marginalizing (\ref{oldenew}) with respect to $(m_1,\dots, m_k)$ yields
$$
p({\bf s}|{\bf n})=p(s_1,\dots, s_{k^*}|{n_1,\dots, n_k})=
$$
\begin{equation}
\label{oldenew2}
= \frac{V_{n+m,k+k^*}}{V_{n,k}} {m \choose m-s}\mathop{\sum_{(m_1,\dots, m_k)}}_{\sum_j m_j=m-s, m_j\geq 0}{m-s \choose{m_1,\dots,m_k}} \prod_{j=1}^k (n_j-\alpha)_{m_j\uparrow}\prod_{j=1}^{k^*}(1-\alpha)_{s_j-1\uparrow}=
\end{equation}
and by (\ref{multi}) 
\begin{equation}
\label{prop1}
=\frac{V_{n+m, k+k^*}}{V_{n,k}}{m \choose m-s}(n-k\alpha)_{m-s\uparrow}\prod_{j=1}^{k^*}(1-\alpha)_{s_j-1\uparrow}.
\end{equation}
%or tutti nei vecchi
 \\\\
{\bf Corollary 6.} Corollary 1, in Lijoi et al. (2008, eq. (10) and (11)),  may be obtained from (\ref{prop1}) by summing over the space of all partitions of $s$ elements in $k^*$ blocks. For $s=\sum_{j=1}^{k^*} s_j$, an application of (\ref{bellalpha}) yields
\begin{equation}
\label{coro1}
Pr(K^*=k^*, S=s|n_1,\dots, n_k)=\frac{V_{n+m, k+k^*}}{V_{n,k}}{m \choose s}(n-k\alpha)_{m-s\uparrow} S_{s,k^*}^{-1,-\alpha}.
\end{equation}
Marginalizing with respect to $k^*$ a probability distribution for the total number $S$ of observations in new blocks is easily derived in terms of generalized Stirling numbers,  
\begin{equation}
\label{obsinnew3}
Pr(S=s|n_1,\dots,n_k)= \frac{1}{V_{n,k}}{m \choose s}(n-k\alpha)_{m-s\uparrow}\sum_{k^*=0}^s {V_{n+m, k+k^*}} S_{s,k^*}^{-1,-\alpha}.
\end{equation}
{\bf Proposition 7.} [Lijoi et. al (2007, eq (4)] { The probability distribution of the number of new blocks $k^*$ may be derived from (\ref{coro1}) and exploiting the convolution formula (\ref{convo})}.\\\\
{\it Proof.} By marginalizing (\ref{coro1}) with respect to $S$, 
\begin{equation}
\label{biom}
Pr(K^*=k^*|n_1,\dots,n_k)=\frac{V_{n+m, k+k^*}}{V_{n,k}}\sum_{s=k^*}^m {m \choose s} (n-k\alpha)_{m-s\uparrow} S_{s,k^*}^{-1,-\alpha}.
\end{equation}
From (\ref{convo}) we know that
\begin{equation}
\label{convo2}
S_{n,k}^{-1, -\alpha, \gamma}= \sum_{s=k}^{n} {n \choose s} S_{s,k}^{-1, -\alpha} (-\gamma)_{n-s \uparrow 1}
\end{equation}
hence
\begin{equation}
\label{29}
Pr(K^*=k^*|n_1,\dots,n_k)=\frac{V_{n+m, k+k^*}}{V_{n,k}}S_{m,k^*}^{-1,-\alpha, -(n-k\alpha)},
\end{equation}
where $S_{m,k^*}^{-1,-\alpha, -(n-k\alpha)}$ is the non-central generalized Stirling number defined in (\ref{noncentralsti}) for scale parameter $\gamma=(n-k\alpha)$. The result in equation (12) in Lijoi et al. (2008), expressed in terms of Stirling numbers, follows from the ratio of (\ref{coro1}) and (\ref{29}) 
$$
Pr(S=s|K^*=k^*, (n_1,\dots, n_k))=\frac{{m \choose s}(n-k\alpha)_{m-s\uparrow} S_{s,k^*}^{-1,-\alpha}}{S_{m,k^*}^{-1,-\alpha, -(n-k\alpha)}}.
$$\\
By interpreting distributions (\ref{obsinnew3}) and (\ref{29}) as posterior distributions in a Bayesian perspective the corresponding expected values give corresponding Bayes estimator under quadratic loss function for $K^*$ and $S$, (see Lijoi et al. 2008, eq. (14) and (15)). \\\\
{\bf Conjecture 8.} In Proposition 2. Lijoi et al. (2008) obtain the following simplification for the expected value of the number $S$ of observations in new blocks, given the placement of the first $n$ customers
$$
E(S|n_1,\dots, n_k)=m\frac{V_{n+1, k+1}}{V_{n,k}}, 
$$
and provide a proof by induction based on properties of conditional expectations under exchangeability. \\
In our setting we conjecture that an alternative proof may be obtained by suitably exploiting the backward recursive equation for the Gibbs weights
$$
V_{n,k}=(n-k\alpha)V_{n+1,k}+V_{n+1,k+1}, 
$$
combined with some already recalled relationships for Stirling numbers.\\\\ 
First notice that by (\ref{bino}), for positive $x$ and $y$ 
$$
\sum_{s=0}^m {m \choose s}\frac{(x)_s (y)_{m-s}}{(x+y)_m}=1.
$$
and by substitution $z=x+1$ and $t=s-1$ it easy to show that
\begin{equation}
\label{iper}
\sum_{s=0}^m s{m \choose s}\frac{(x)_s (y)_{m-s}}{(x+y)_m}=m\left(\frac{x}{x+y}\right)\sum_{t=0}^{m-1} {m-1 \choose t} \frac{(z)_t (y)_{m-1-t}}{(z+y)_{m-1}}=m\left(\frac{x}{x+y}\right).
\end{equation}
Therefore, if one is able to show that the probability distribution of $S$  
\begin{equation}
\label{uno}
Pr(S=s|n_1,\dots,n_k)=\frac{1}{V_{n,k}}{m \choose s}(n-k\alpha)_{m-s\uparrow}\sum_{k^*=0}^s {V_{n+m, k+k^*}} S_{s,k^*}^{-1,-\alpha}
\end{equation}
can be reduced to 
\begin{equation}
\label{tre}
Pr(S=s|n_1,\dots, n_k)={m \choose s} \frac{((n-k\alpha)V_{n+1,k})_{m-s\uparrow}(V_{n+1,k+1})_{s\uparrow}}{(V_{n,k})_{m_\uparrow}}
\end{equation}
the result would follow. We conjecture that the intermediate step would be to prove that the backward recursive equation sufficies to show that (\ref{uno}) equates
\begin{equation}
\label{due}
{m \choose s}(n-k\alpha)_{m-s\uparrow}\frac{V_{n+m,k}}{V_{n,k}}\sum_{k^*=0}^{s} (V_{n+1, k+1})_{k^* \uparrow \alpha} S_{s,k^*}^{-1,-\alpha},
\end{equation}
in fact, by means of (\ref{unoalpha}), it would reduce to (\ref{tre}).
We remark that, up to now, this is a conjecture and equivalence of (\ref{uno}), (\ref{due}) and (\ref{tre}) is still to be proved.\\\\
{\bf Example 9}. The previous conjecture is inspired by the particular case of the two parameter Poisson-Dirichlet model (Pitman and Yor, 1997). It is known that it belongs to the Gibbs class of type $\alpha \in (0,1)$ since arises by mixing the stable conditioned $PK(\rho_\alpha|t)$ model by $\gamma(t)=\frac{\Gamma(\theta+1)}{\Gamma(\theta/\alpha +1)}t^{-\theta}f_{\alpha}(t)$. The Gibbs weights are known to be
$$V_{n,k}=\frac{(\theta +\alpha)_{k-1\uparrow \alpha}}{(\theta +1)_{n-1\uparrow}}.$$
By exploiting the recursive equation some algebra shows that
$$
\frac{V_{n+1,k+1}}{V_{n,k}}=1-\frac{(n-k\alpha)V_{n+1,k}}{V_{n,k}}=\frac{(\theta +k\alpha)}{(\theta +n)}
$$
and 
$$
\frac{V_{n+m,k}}{V_{n,k}}=\frac{1}{(\theta +n)_{m}}.
$$
By (\ref{multiplicative}) and (\ref{unoalpha}), the sum in (\ref{obsinnew3}) results
$$
\frac{1}{(\theta +1)_{n+m-1}}\sum_{k^*=0}^s (\theta +\alpha)_{k+k^*-1\uparrow \alpha} S_{s,k^*}^{-1, -\alpha}=
\frac{(\theta +\alpha)_{k-1 \uparrow \alpha}}{(\theta +1)_{n+m-1}}\sum_{k^*=0}^s  (\theta + k \alpha)_{k^* \uparrow \alpha} S_{s, k^*}^{-1, -\alpha}= \frac{(\theta +\alpha)_{k-1 \uparrow \alpha}}{(\theta +1)_{n+m-1}} (\theta +k\alpha)_{s \uparrow}.
$$
hence equation (\ref{obsinnew3})  specializes as\\\
$$
Pr(S=s| n_1,\dots,n_k)
={m \choose s}(n -k \alpha)_{m-s}\uparrow \frac{V_{n+m,k}}{V_{n,k}}(\theta +k\alpha)_s\uparrow=
 {m \choose s}\frac{(n-k\alpha)_{m-s\uparrow}(\theta +k\alpha)_{s \uparrow}}{(\theta +n)_{m \uparrow}}
$$
and by (\ref{iper}) the expected value results
$$E(S|n_1,\dots, n_k)=m\frac{(\theta +k\alpha)}{\theta +n}.
$$
{\bf Corollary 10.} [Proposition 4, Lijoi et al. (2008)] The probability that the $m$ new customers don't seat at a subset of $(k-r)$ old tables arises from (\ref{oldenew}) by summing over the ways to choose $s$ customers from the $m$ of the new group, by summing over the ways to partition $s$ customers in a subset of $k^*$ new tables, for $k^* \geq 1$, and over the ways to allocate $m-s$ customers in at most $r$ old tables. Notice that the operations of partitioning and allocating differ for the fact that the  blocks of partitions cannot be empty while allocation in a fixed number of blocks can result in a certain number of blocks remaining empty. From (\ref{oldenew}) we obtain

\begin{equation}
\label{backward}
\sum_{k^*=1}^m \frac{V_{n+m, k+k^*}}{V_{n,k}}\frac{1}{k^*!} \sum_{s=k^*}^m {m \choose s} \sum_{s_1,\dots, s_{k^*}, \sum_j s_j=s} {s \choose {s_1 \cdots s_{k^*}}} \prod_{j=1}^{k^*} (1 -\alpha)_{s_j -1 \uparrow} \times
\end{equation}
\begin{equation}
\times \mathop{\sum_{(m_1,\dots, m_r),\sum_j m_j=m-s}} {m-s \choose {m_1 \cdots m_r}} \prod_{j=1}^r (n_j -\alpha)_{m_j \uparrow}=
\end{equation}
which simplifies to  
$$
\sum_{k^*=1}^m \frac{V_{n+m, k+k^*}}{V_{n,k}}\sum_{s=k^*}^m {m \choose s} S_{s,k^*}^{-1,- \alpha} (\sum_{j=1}^{r} n_j -r\alpha)_{m-s\uparrow}=
$$
and the definition of non central generalized Stirling numbers (cfr. eq. (\ref{convo})) yields 
$$
=\sum_{k^*=1}^m \frac{V_{n+m, k+k^*}}{V_{n,k}} S_{m, k^*}^{-1,-\alpha,(r\alpha -\sum_j n_j)}.
$$

\subsection {The reproducibility of the Gibbs structure from "deletion of classes" property of PK models}

In Section 3.1 Lijoi et al. (2008) point out a "reproducibility" property of the exchangeable Gibbs partitions of type $\alpha$ which motivates the definition of {\it conditional Gibbs structures}. Here we investigate the relationship between this property and the {\it deletion of classes} property of Poisson-Kingman models introduced in Pitman (2003).\\\\
{\bf Definition 11.} [Deletion of classes, Pitman (2003)] Given a random partition $\Pi$ of $\mathbb{N}$, the operator {\it deletion of the first $k$ classes} is as follows: First let $\Pi_k^*$ be the restriction of $\Pi$ to $H_k:=\mathbb{N}-G_1-\cdots-G_k$ where $G_1,\dots, G_k$ are the first $k$ classes of $\Pi$ in order of their least elements, then derive $\Pi_k$ on $\mathbb{N}$ from $\Pi^*_k$ on $H_k$ by renumbering the points of $H_k$ in increasing order. \\\\
{\bf Proposition 12.} The {conditional EPPF in Definition 2. of Lijoi, et al. (2008) yields the EPPF of a random partition obtained by the operation of {\it deletion of classes}}.\\\\  
{\it Proof.} From our group sequential construction the EPPF of $\Pi_k$  is given by 
\begin{equation}
\label{dele}
p({\bf s| m, n})=\frac{p({\bf m,s|n})}{p({\bf m|n})}.
\end{equation}
Recall from equation (20) that $s < m$ customers sit at $k^*$ {\it new} tables in configuration $(s_1,\dots,s_{k^*})$ and the remaining $m-s$ customers sit at the {\it old} tables in configuration $(m_1,\dots, m_k)$ for $\sum_{j=1}^{k} m_j= m-s$, $1 \leq s \leq m$, $\sum_{j=1}^{k^*} s_j=s$, $m_j \geq 0$, $s_j \geq 1$ with probability
\begin{equation}
\label{oldnew}
p_{{\bf m, \bf s}}({\bf n})=p({\bf m}, {\bf s}| {\bf n})=\frac{p(n_1+m_1, \dots, n_k+m_k, s_1, \dots, s_{k^*})}{p(n_1,\dots, n_k)},\\
\end{equation}
and for Gibbs partitions of type $\alpha$ this yields (cfr. eq. (\ref{oldenew})
\begin{equation}
\frac{V_{n+m,k+k^*}}{V_{n,k}}\prod_{j=1}^k (n_j-\alpha)_{m_j\uparrow}\prod_{j=1}^{k^*}(1-\alpha)_{s_j-1\uparrow}.
\end{equation}
The denominator in (\ref{dele}) is obtained by marginalizing (\ref{oldnew}) with respect to ${\bf s}=(s_1,\dots, s_{k^*})$, i.e. by summing over all the ways to partition $s$ observations in $k^*$ new tables for every $k^*$, i.e.
%\sum_{k^*=1}^s \sum_{(s_1,\dots, s_{k^*})} p({\bf s,m|n})=\prod_{j=1}^{k} (n_j-\alpha)_{m_j\uparrow} \sum_{k^*=1}^s  \frac{V_{n+m,k+k^*}}{V_{n,k}} \sum_{(s_1,\dots, s_{k^*})} \prod_{j=1}^{k^*} (1-\alpha)_{s_j-1\uparrow}
$$
p({\bf m|n})=\prod_{j=1}^{k} (n_j-\alpha)_{m_j\uparrow} \sum_{k^*=1}^s  \frac{V_{n+m,k+k^*}}{V_{n,k}} S_{s,k^*}^{-1,\-\alpha},
$$
hence the EPPF of $\Pi_k$ is given by  
%obtained by deletion of classes is given by
$$
p({\bf s|m,n})=\frac{V_{n+m,k+k^*}}{\sum_{k^*=1}^sV_{n+m,k+k^*}S_{s,k^*}^{-1,-\alpha} }\prod_{j=1}^{k^*} (1-\alpha)_{s_j-1\uparrow}
$$
which agrees with the result in Proposition 3. in Lijoi et al. (2008).\\\\ 
{\bf Remark 13.} Notice that in Lijoi et al. (2008) the result is obtained conditioning on the number $s$ of customers in new tables. Given the size $m$ of the {\it new} group, this is equivalent to conditioning to the number $m-s$ of the observations in old blocks, i.e. to the vector $(m_1,\dots, m_k)$ as in our result.\\\\\
Now from Proposition 7. in Pitman (2003) if $\Pi$ is a Poisson Kingman $PK(\rho, \gamma)$ partition of $\mathbb{N}$, and $\Pi_k$ is derived from $\Pi$ by deletion of the first $k$ classes, then  $\Pi_k$ is a $PK(\rho, \gamma_k)$ of $\mathbb{N}$ where $\gamma_k$ is given by $\gamma_k=\gamma Q^k$, where $Q$ is the Markov transition operator on $(0,\infty)$ 
$$
Q(t,dv)=\rho(t-v)(t-v)t^{-1}f(v)1\{0<v<t\}dv
$$ 
and $f(t)$ is the probability density corresponding the Levy density $\rho$. Since from  Gnedin and Pitman (2006) (cfr. item (iii) of Th. 2), we know that an exchangeable partition belongs to the class of Gibbs form of type $\alpha \in (0,1)$ if and only if is a mixture of Poisson-Kingman models 
for some mixing density $\gamma$, and since $\Pi_k$ obtained by deletion of classes of $PK(\rho_\alpha, \gamma)$ produces $PK(\rho_\alpha, \gamma_k)$ which is still of Gibbs form, it should follows that the reproducibility of the Gibbs class also holds for infinite conditional structures. This point, which seems to contradict Corollary 2. in Lijoi et al. (2008), deserves futher investigation that we postpone to a future paper. \\\\

\section*{References}
\newcommand{\bibu}{\item \hskip-1.0cm}
\begin{list}{\ }{\setlength\leftmargin{1.0cm}}

%\bib \textsc{Aalen, O. O.} (1992) Modelling heterogeneity in survival analysis by the compound Poisson distribution. {\it Ann. Appl. Probab.} 2,  951-972.

%\bib \textsc{Antoniak, C. E.} (1974) Mixtures of Dirichlet processes with applications to Bayesian nonparametric problems. {\it Ann. Statist.} 2, 1152-1174. 

%\bib \textsc{Barndorff-Nielsen, O. E. and Shepard, N.} (2001)  Normal modified stable processes.  {\it Th. Probab. Math. Statist.}, 65, 1-19.

%\bib \textsc{Brix, A.} (1999) Generalized Gamma measures and shot-noise Cox processes. {\it Adv. Appl. Probab.}, 31, 929--953.

%\bibu \textsc{Cerquetti, A.} (2007) A note on Bayesian nonparametric priors derived from exponentially tilted Poisson-Kingman models. {\it Stat. \& Prob. Lett} 77, 1705-1711.

%\bibu \textsc {Cerquetti, A.} (2007) On a Gibbs characterization of normalized generalized Gamma processes. arXiv 0707:3408. 

\bibu  \textsc {Charalambides, C. A.} (2005) {\it Combinatorial Methods in Discrete Distributions}. Wiley, Hoboken NJ.

%\bibu \textsc{Engen, S.} (1978) {\it Stochastic abundance models}. Chapman \& Hall, London.

%\bibu \textsc{Ewens, W. and Tavar\'e S.} (1995) The Ewens sampling formula. In Multivariate discrete distributions (Johnson, N.S., Kotz, S. and Balakrishnan, N. eds.). Wiley, NY.

%\bibu \textsc{Ferguson, T. S.} (1973)  A Bayesian analysis of some nonparametric problems. {\it Ann. Statist.}, 1, 209--230.

%\bib \textsc{Fisher, R.A., Corbet, A.S. and Williams, C.B.} (1943) The relation beteween the number of species and the number of individuals in a random sample of an animal population. {\it J. Animal. Ecol.}, 12, 42--58.

\bibu \textsc{Gnedin, A. and Pitman, J. } (2006) {Exchangeable Gibbs partitions  and Stirling triangles.} {\it Journal of Mathematical Sciences}, 138, 3, 5674--5685. 

%\bibu \textsc{Hansen, B. and Pitman, J.} (2000) Prediction rules for exchangeable sequences related to species sampling. {\it Statistics \& Probability Letters}, 46, 251--256.

\bibu \textsc{Ho, M-W, James, L.F. and Lau, J.W.} (2007) Gibbs partitions (EPPF's) derived from a stable subordinator are Fox H - And Meijer G - Transforms. arXiv:0708.0619v2 [math.PR]

\bibu \textsc{Hsu, L. C, \& Shiue, P. J.} (1998) A unified approach to generalized Stirling numbers. {\it Adv. Appl. Math.}, 20, 366-384.

%\bib \textsc{Ho, M., James, L.F. and Lau, J.W} (2006) Coagulation Fragmentation laws induced by general coagulations of two-parameter Poisson-Dirichlet processes. {\it arXiv:math.PR/0601608}

%\bib \textsc{Hougaard, P} (1986) Survival models for hetereneous populations derived from stable distributions. {\it Biometrika}, 73, 387--396.

\bibu \textsc{Ishwaran, H. \& James, L. F.} (2003) Generalized weighted Chinese restaurant processes for species sampling mixture models. {\it Statist. Sinica}, {13}, 1211--1235.

%\bibu \textsc{Ishwaran, H. \& Zarepour, M.} (2002) Exact and approximate sum-represen\-tations for the Dirichlet process. {\it Can. J. Statist.} 30, 269-283. 

%\bibu \textsc{James, L. F.} (2002). Poisson process partition calculus with applications to exchangeable models and Bayesian Nonparametrics. {\it arXiv:math.ST/0205093}.

%\bib \textsc{James, L. F.} (2003) A simple proof of the almost sure discreteness of a class of random measures. {\it Statist. \& Probab. Lett.}, 65, 363-368.

%\bib \textsc{James, L.F., Lijoi, A. and Pr\"unster I.} (2005) Bayesian inference via classes of normalized random measures. {\it arXiv:math.ST/0503394}.

\bibu \textsc{Kerov, S.} (1995) Coherent random allocations and the Ewens-Pitman sampling formula. PDMI Preprint, Steklov Math. Institute, St. Petersburg. 

%\bibu \textsc{Kingman, J.F.C} (1967) Completely random measures. {\it Pacific J. Math.}, 21, 59-78

\bibu \textsc{Kingman, J.F.C.} (1975) Random discrete distributions. {\it J. Roy. Statist. Soc. B}, 37, 1--22. 

\bibu \textsc{Kingman, J.F.C} (1978) The representation of partition structure.  {\it J. London Math. Soc.} 2, 374--380.

%\bibu \textsc{Lijoi, A., Mena, R. and Pr\"unster, I.} (2005) Hierarchical mixture modeling with normalized Inverse-Gaussian priors. {\it JASA}, vol. 100, 1278--1291.

\bibu \textsc{Lijoi, A., Mena, R. and Pr\"unster, I.} (2007) Bayesian nonparametric estimation of the probability of discovering new species  {\it Biometrika}, 94, 769--786.

\bibu \textsc{Lijoi, A., Pr\"unster, I. and Walker, S.G.} (2008) Bayesian nonparametric estimator derived from conditional Gibbs structures. {\it Annals of Applied Probability}, (To appear)

\bibu \textsc {Normand, J.M.} (2004) Calculation of some determinants using the $s$-shifted factorial. {\it J. Phys. A: Math. Gen.} 37, 5737-5762.

\bibu \textsc{Perman, M., Pitman, J, \& Yor, M.} (1992) Size-biased sampling of Poisson point processes and excursions. {\it Probab. Th. Rel. Fields}, 92, 21--39.

%\bibu \textsc{Pitman, J.} (1996) Some developments of the Blackwell-MacQueen urn scheme. In T.S. Ferguson, Shapley L.S., and MacQueen J.B., editors, {\it Statistics, Probability and Game Theory}, volume 30 of {\it IMS Lecture Notes-Monograph Series}, pages 245--267. Institute of Mathematical Statistics, Hayward, CA.

\bibu \textsc{Pitman, J.} (2003) {Poisson-Kingman partitions}. In D.R. Goldstein, editor, {\it Science and Statistics: A Festschrift for Terry Speed}, volume 40 of Lecture Notes-Monograph Series, pages 1--34. Institute of Mathematical Statistics, Hayward, California.

\bibu \textsc{Pitman, J.} (2006) {\it Combinatorial Stochastic Processes}. Ecole d'Et\'e de Probabilit\'e de Saint-Flour XXXII - 2002. Lecture Notes in Mathematics N. 1875, Springer.

\bibu \textsc{Pitman, J. and Yor, M.} (1997) The two-parameter Poisson-Dirichlet distribution derived from a stable subordinator. {\it Ann. Probab.}, 25:855--900.

%\bib \textsc{Regazzini, E., Lijoi, A. and Pr\"unster, I.} (2003) Distributional results for means of random measures with independent increments. {\it Ann. Statist.}, 31, 560--585.

%\bib \textsc{Sato, K.} (1999) {\it L\'evy processes and infinitely divisible distributions}. Cambridge University Press.

%\bib \textsc{Seshadri, V. (1993)} {\it The inverse Gaussian distribution}. Oxford University Press, New York. 

%\bib \textsc{Sibisi, S. \& Skilling, J.} (1997) Prior distributions on measure space. {\it J.R. Statist. Soc.} B, 59, 1, 217--235.

%\bib \textsc{Steutel F, W, \& Van Harn, K.} (2004) {\it Infinite divisibility of probability distributions on the real line}. {Monographs and Textbooks in Pure and Applied Mathematics}, 259. Marcel Dekker, New York.

\bibu \textsc{Toscano, L.} (1939) Numeri di Stirling generalizzati operatori differenziali e polinomi ipergeometrici. {\it Comm. Pontificia Academica Scient. } 3:721-757.

\end{list}

\end{document}